\documentclass[]{article}

\usepackage{fullpage}
\usepackage[hidelinks]{hyperref} 

\usepackage{amsmath,amsthm,amssymb}
\usepackage{mathtools}
\usepackage{enumitem}

\usepackage{subcaption} 

\numberwithin{equation}{section}

\newtheorem{theorem}{Teorema}[section]
\theoremstyle{definition}
\newtheorem{remark}[theorem]{Observación}

\usepackage{xifthen}
\newcommand{\dC}[1][]{\prescript{C}{}{\!D^{\ifthenelse{\isempty{#1}}{\alpha}{#1}}}}
\newcommand{\dRL}[1][]{\prescript{RL}{}{\!D^{\ifthenelse{\isempty{#1}}{\alpha}{#1}}}}
\newcommand{\dGL}[1][]{\prescript{GL}{}{\!D^{\ifthenelse{\isempty{#1}}{\alpha}{#1}}}}
\newcommand{\Laplace}{{\mathcal L}}

\title{Fractional Derivatives: \\ an extension of classical analysis to non-integer orders}

\author{Félix del Teso%
\thanks{Departamento de Matemáticas, Universidad Autónoma de Madrid. \texttt{felix.delteso@uam.es}}
\thanks{ICMAT-Instituto de Ciencias Matemáticas, CSIC-UAM-UC3M-UCM
} 
\and 
David Gómez-Castro%
\thanks{Departamento de Matemáticas, Universidad Autónoma de Madrid. \texttt{david.gomezcastro@uam.es}} 
\thanks{ICMAT-Instituto de Ciencias Matemáticas, CSIC-UAM-UC3M-UCM
}}

\begin{document}

\maketitle

\begin{abstract}
    This article provides an accessible introduction to fractional derivatives, a concept that extends classical calculus by allowing derivatives of non-integer order. It explores both the fundamental definitions and some of the most relevant properties and applications of this mathematical tool.
    It was originally published in Spanish in the \emph{Gaceta de la Real Sociedad Española} as \cite{delTesoGomezCastro2025Gaceta} and automatically translated using Github copilot.
\end{abstract}

\section{Introduction}\label{sec:intro}

The origin of fractional differentiation and integration operators dates back to the 17th century. In 1675, G. Leibniz introduced the concept of the $n$-th derivative (for $n\in \mathbb{N}$) of a function. Later, in 1695, he exchanged correspondence with the Marquis de L'Hôpital (see \cite{leibniz1962}), where the latter raised the question of the possible meaning of differentiating $n = 1/2$ times, to which Leibniz replied:
\begin{quote}
    \textit{Il y a de l'apparence qu'on tirera un jour des consequence bien utiles de ces paradoxes, car il n'y a gueres de paradoxes sans utilité.}
\end{quote}
which in English reads:
\begin{quote}
    There is an appearance that one day very useful consequences will be drawn from these paradoxes, for there are hardly any paradoxes without utility.
\end{quote}
In one of the letters, Leibniz proposed that a possible definition of the fractional derivative $\frac{d^\alpha}{dx^\alpha}$ of order $\alpha\in(0,1)$ of the function $f(x)=x$ should satisfy
\begin{equation*}
    \frac{d^\alpha f}{d x^\alpha }(x)= x^{1-\alpha},
\end{equation*}
which in the particular case $\alpha=1/2$ would give $\frac{d^{1/2} }{d x^{1/2} }f(x)= x^{1/2}$ (as we will see, it is reasonable to think that some constants are missing in the formulas above).
The first formal reference to a concept of fractional derivative appears years later in the article \cite{euler1738} by L. Euler,
who defines it through his formula for interpolating factorial numbers between positive integers. In this article we will not give a detailed historical introduction beyond this brief note. For further reading in this direction we recommend \cite{MorilloRiera-Segura2022BreveIntroduccionCalculo}.

This introductory text on fractional derivatives is organized as follows. In Section \ref{sec:gamma} we briefly recall the definition and some properties of the Gamma function $\Gamma$, which generalizes the factorial concept to non-integer orders. After that, in Section \ref{sec:primera def} we provide a first informal definition based on the correspondence between L'Hôpital and Leibniz, applied to monomials, as a possible approach to the fractional derivative. Section \ref{sec:def formales} is devoted to defining the fractional (integration and differentiation) operators of Riemann-Liouville, Caputo, and Grünwald-Letnikov. Moreover, we will see how they are related to each other and present some of their most relevant properties. The following sections, devoted to various topics related to fractional derivatives, are accompanied by the associated classical theory that is typically studied in a first course on ordinary differential equations in a science degree.
Section \ref{sec: teoremas fundamentales} focuses on deriving the fundamental theorems of fractional calculus, that is, how fractional derivatives and integrals interact. Next, in Section \ref{sec:ODE}, we present and study fractional differential equations of various types. As expected, most of these equations cannot be solved explicitly, which leads us in Section \ref{sec:metodos numericos} to introduce some numerical methods. In Section \ref{sec:sistemas} we make some comments about systems of fractional differential equations. Finally, in Section \ref{sec:apli} we give examples of applications. In Section \ref{sec:tautocrona} we describe the tautochrone problem, whose solution is given by the solution of the fractional equation
\begin{align*}
    \frac{d^{1/2}s}{dy^{1/2}} (y) & = T(y) \sqrt{\frac{2g}{\pi}} , \quad \text{for } y>0 \\
    s(0)                          & =0,
\end{align*}
where $T(y)$ is the chosen fall time (constant in the case of the tautochrone), and $s$ is the arc-length function.
We also present, in Section \ref{sec:finanzas}, a more involved application, perhaps more important nowadays, to asset pricing models.

\subsection{The \texorpdfstring{$\Gamma$}{Gamma} function}
\label{sec:gamma}
\begin{figure}
    \centering
    \includegraphics[width=.6\textwidth]{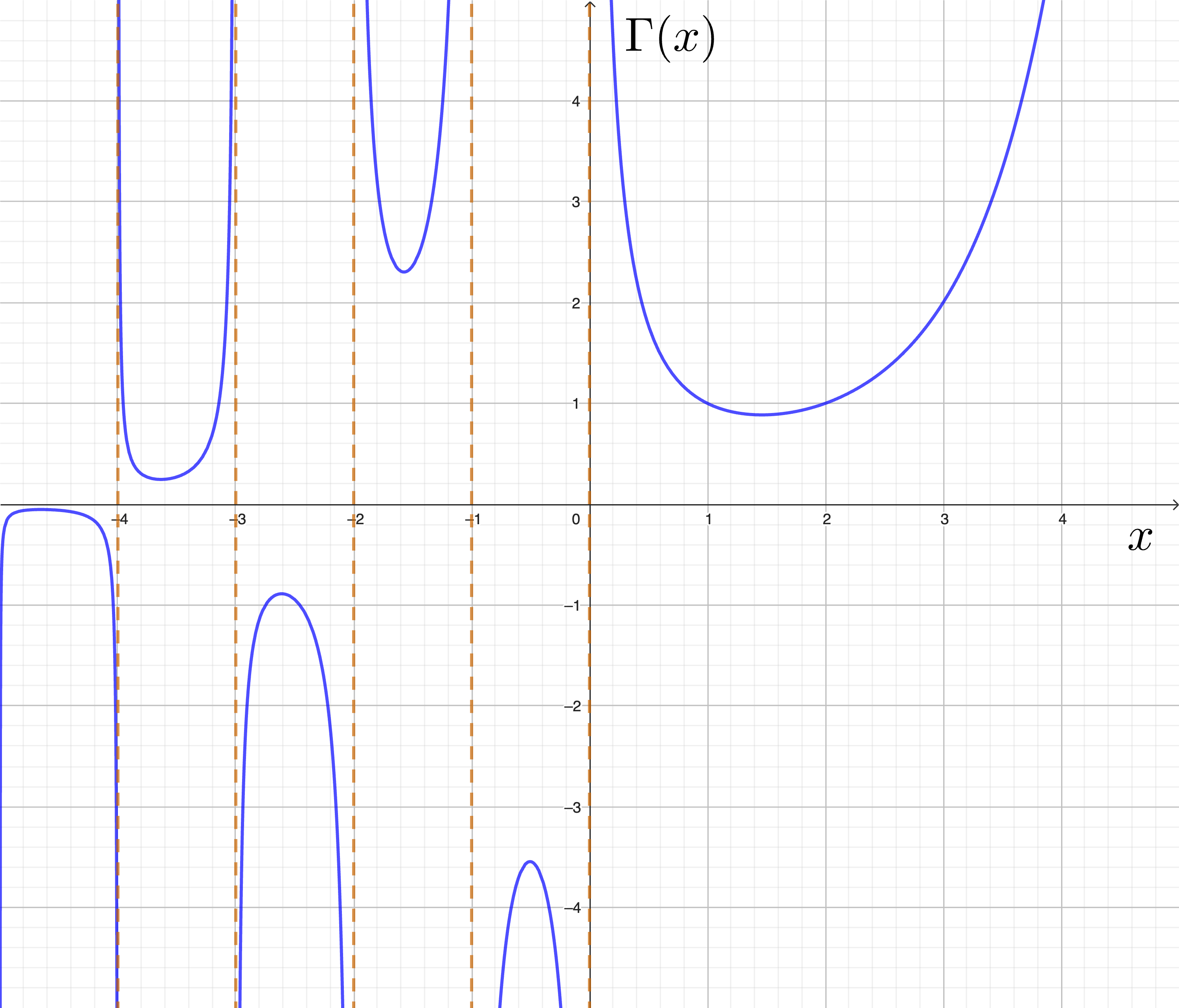}
    \caption{Graph of the Gamma function $\Gamma(x)$.}
    \label{figura:gamma}
\end{figure}

Since it will be useful throughout the exposition, we now give a brief overview of Euler’s famous Gamma function. It satisfies $\Gamma : (0,\infty) \to (0,\infty)$, $\Gamma(n) = (n-1)!$ for $n\in \mathbb{N}$, and moreover it is smooth and increasing for $z \ge 1$. Its exact expression is
\begin{equation*}
    \Gamma(z) \coloneqq \int_0^\infty s^{z-1} e^{-s} ds, \qquad \text{for all } z > 0.
\end{equation*}
In what follows we will use the following known properties:
\begin{enumerate}
    \item Iteration property: $\Gamma(z) z = \Gamma(z + 1)$.
    \item Note that, for $z > 0$, we have
          \begin{equation}
              \label{eq:gamma(0)}
              \Gamma(z) \ge \int_0^1 \sigma^{z-1} e^{-1} d\sigma = \frac{e^{-1}}{z} \to +\infty, \qquad \text{cuando } z \to 0^+.
          \end{equation}
\end{enumerate}
For $z \in (-\infty,0) \setminus \mathbb Z$ we can define $\Gamma$ from the iteration property. For example, we define
\[
    \Gamma(z) \coloneqq \frac{\Gamma(1+z)}{z} \quad \text{ for } z \in (-1,0).
\]
Note that $|\Gamma(-n^\pm)| = \infty$ if $n \in \mathbb N$. This divergent character can be clearly observed in Figure \ref{figura:gamma}.

We will also need the Beta function, defined for $z_1,z_2>0$ as
\begin{equation*}
    \mathrm B(z_1, z_2) = \int_0^1 t^{z_1-1} (1 - t)^{z_2 -1} dt.
\end{equation*}
It enjoys the following property relating it to the Gamma function:
\begin{equation}\label{eq:relgammabeta}\mathrm B(z_1, z_2) = \frac{\Gamma(z_1)\Gamma(z_2)}{\Gamma(z_1 + z_2)}.
\end{equation}
\subsection{A first informal definition based on the correspondence between L'H\^{o}pital and Leibniz}
\label{sec:primera def}

Based on Leibniz’s idea (explained in Section \ref{sec:intro}) of how a fractional derivative should act on monomials $f_\beta(t)=t^{\beta}$ for $\beta\not=0$, we look for an operator $D^\alpha$ with $\alpha\geq0$ such that
\begin{align}\label{eq:propgrado}
    D^\alpha [f_\beta] (t)=  C(\alpha,\beta)f_{\beta-\alpha} (t),
\end{align}
for a certain family of constants $C=C(\alpha,\beta)$. That is, “differentiating $\alpha$ times” a monomial of degree $\beta$ reduces the degree to $\beta-\alpha$ (up to a multiplicative constant). Property \eqref{eq:propgrado}, by itself, is too unrestrictive and allows too many possible definitions (depending on the choice of the family of constants). Since we would like these derivatives to be a generalization of classical derivatives (i.e. when $\alpha\in \mathbb{N}$), it is natural to require the composition (semigroup) property: given $\alpha_1,\alpha_2\geq0$, we would like
\begin{align}
    \label{eq:comp}
    D^{\alpha_{2}} [D^{\alpha_{1}}[f_\beta]]=D^{\alpha_{1}+\alpha_{2}}[f_\beta].
\end{align}
On the one hand,
\[
    D^{\alpha_{2}} [D^{\alpha_{1}}[f_\beta]]= C(\alpha_1,\beta)D^{\alpha_{2}}[f_{\beta-\alpha_1}]=C(\alpha_1,\beta)C(\alpha_2,\beta-\alpha_1) f_{\beta-\alpha_1-\alpha_2},
\]
and on the other,
\[
    D^{\alpha_{1}+\alpha_{2}}[f_\beta]= C(\alpha_1+\alpha_2,\beta)f_{\beta-\alpha_1-\alpha_2}.
\]
Therefore, property \eqref{eq:comp} requires that our choice of constants satisfy
\begin{align}\label{eq:propconstantes}
    C(\alpha_1+\alpha_2,\beta)=C(\alpha_1,\beta)C(\alpha_2,\beta-\alpha_1).
\end{align}
The existence of such constants (not completely determined by this condition) is easy to verify. For this we can use as a hint the value of the constant in the case $\alpha=n\in \mathbb{N}_{>0}$ and $\beta  \in \mathbb{N}_{>0}$, that is, when we work with classical derivatives. In this case we know that
\[
    D^{n}[f_{\beta}]=\beta\cdot(\beta-1)\cdot \cdots \cdot (\beta-n+1) f_{\beta-n}= \frac{\beta!}{(\beta-n)!}f_{\beta-n} =  \frac{\Gamma(\beta+1)}{\Gamma(\beta+1-n)}f_{\beta-n},
\]
where we have used the definition of the Gamma function described in the previous subsection.
Therefore, generalizing for any $\alpha\geq0$ and $\beta \in \mathbb{R} \setminus \{-1,-2,\cdots\}$, it seems natural to choose
\[
    C(\alpha,\beta):= \frac{\Gamma(\beta+1)}{\Gamma(\beta+1-\alpha)}.
\]
The restriction on $\beta$ is due to the fact that the Gamma function is not defined (because of its divergent nature) at nonpositive integers (see again the previous subsection), so we avoid them in the numerator (they can appear in the denominator, making $C(\alpha,\beta)=0$ in that case).
Let us check that this choice indeed satisfies \eqref{eq:propconstantes}:\footnote{The same argument shows that, for any function $F$, the family of constants $C(\alpha,\beta):= F(\beta)/F(\beta-\alpha)$ (with restrictions on $\alpha, \beta$ depending on the domain of $F$) also satisfies \eqref{eq:propconstantes}. Nevertheless, choosing $F(\cdot)=\Gamma(\cdot+1)$ is the natural and most interesting option.}
\begin{align*}
    C(\alpha_1,\beta)C(\alpha_2,\beta-\alpha_1) & = \frac{\Gamma(\beta+1)}{\Gamma(\beta+1-\alpha_1)} \cdot \frac{\Gamma(\beta-\alpha_1+1)}{\Gamma(\beta+1-\alpha_1+\alpha_2)} \\
                                                & =\frac{\Gamma(\beta+1)}{\Gamma(\beta+1-\alpha_1+\alpha_2)}                                                                  \\
                                                & = C(\alpha_1+\alpha_2,\beta).
\end{align*}

\begin{remark}
    Very regular functions can be expressed as power series. As we have seen, it is easy to define the fractional derivative of monomials and then attempt to express the fractional derivative of more general functions. But this approach, while natural, is not the most effective.
\end{remark}

\section{Definitions of fractional operators}\label{sec:def formales}

We will be interested in two types of fractional operators: derivatives and integrals. We start with the latter, since they arise more naturally.

\subsection{Riemann–Liouville fractional integral}\label{sec:intfrac}

We present two equivalent definitions of the fractional integral that appeared at different historical moments.
Both are certainly well defined for functions $f \in C([0,\infty))$, but they can also be defined for functions with less regularity, for example $C((0,\infty))$ and bounded. Unless this could lead to confusion, we will not specify the regularity of the functions.

\paragraph{Liouville and Riemann's proposal.}
This formalization, based on a generalization of the integral (or antiderivative) to arbitrary orders, was introduced by J. Liouville in \cite{liouville1832memoire}, and by B. Riemann in a manuscript published posthumously as part of his collected works \cite{riemann1898oeuvres}.

Let $I$ be the integral (primitive) operator
\[
    If (t) \coloneqq \int_0^t f(s) ds.
\]

For reasons that will be immediately clear, set $I_1 \coloneqq I$.
When we want to compute $I_2 \coloneqq I \circ I$ we can use Fubini's Theorem to write
\begin{align*}
    I_2f (t) = \int_0^t I_1f (s) ds = \int_0^t \int_0^s f(\sigma) d\sigma ds = \int_0^t f(\sigma) \int_\sigma^t  ds d \sigma = \int_0^t f(\sigma) (t-\sigma) d \sigma.
\end{align*}
Inductively, it is easy to prove that $I_n \coloneqq I \circ \overset{n} \cdots \circ I$ can be written as
\[
    I_n f (t) = \int_0^t f(\sigma) K_n(t-\sigma) d \sigma, \quad \text{where} \quad K_n(t) \coloneqq \frac{t^{n-1}}{(n-1)!}=\frac{t^{n-1}}{\Gamma(n)},
\]
where we used the convention $0!=1$.
This gives the intuition, and so it happened historically, that for every $\alpha>0$ one can define the \textit{Riemann–Liouville fractional integral} of order $\alpha>0$ as
\begin{align}\label{eq:defintfrac}
    I_\alpha f (t) \coloneqq \int_0^t f(\sigma) K_\alpha(t-\sigma) d \sigma, \quad \text{where} \quad  K_\alpha (t) \coloneqq \frac{t^{\alpha -1}}{\Gamma(\alpha)}.
\end{align}

Let us see that, indeed, this definition of fractional integral behaves as expected when acting on monomials $f_\beta(t)=t^{\beta}$ with $\beta \geq0$. Using the change of variables $\sigma=ts$, we have
\begin{equation*}
    \begin{aligned}
        I_\alpha f_\beta(t) & = \frac{1}{\Gamma(\alpha)} \int_0^t \sigma^\beta (t-\sigma)^{\alpha - 1} d\sigma = \frac{t^{\beta+\alpha}}{\Gamma(\alpha)} \int_0^1 s^\beta (1-s)^{\alpha - 1}ds  = \frac{\mathrm{B}(\beta+1, \alpha)}{\Gamma(\alpha)} t^{\beta + \alpha},
    \end{aligned}
\end{equation*}
where $\mathrm{B}$ is the Beta function introduced in subsection \ref{sec:gamma}.
Thus, given property \eqref{eq:relgammabeta} of the function $\mathrm{B}$, we obtain
\begin{equation}
    \label{eq:Ialpha de monomio}
    I_\alpha f_\beta(t) = \frac{\Gamma(\beta + 1)}{\Gamma(\alpha + \beta + 1)} t^{\beta + \alpha}.
\end{equation}

We note that integrating $\alpha$ times a monomial of degree $\beta$ increases its degree to $\beta+\alpha$, which generalizes what is known in classical calculus.

Clearly, for $\alpha\in \mathbb N_{>0}$ the Riemann–Liouville integral coincides with the classical integral. However, for $\alpha=0$, the divergent character of the Gamma function near zero (see Figure \ref{figura:gamma}) does not allow defining it directly. We are going to check that, as desirable, we have
\begin{equation}\label{def:I0}
    I_0f(t)\coloneqq\lim_{\alpha\to0^+} I_\alpha f(t)=f(t).
\end{equation}
Given a differentiable function $f$, using the expansion $f(\sigma)=f(t)+f'(\xi_{\sigma,t})(\sigma-t)$, where $\xi_{\sigma,t}$ takes a value between $\sigma$ and $t$, and the iteration property of $\Gamma$, we have
\begin{align*}
    I_\alpha f (t) & = \frac{f(t)}{\Gamma(\alpha)}\int_0^t (t-\sigma)^{\alpha-1}d \sigma- \frac{1}{\Gamma(\alpha)} \int_0^t f'(\xi_{\sigma,t}) (t-\sigma)^\alpha d \sigma=\frac{f(t)}{\Gamma(\alpha+1)}+R_\alpha(t).
\end{align*}
Clearly,
$
    \lim_{\alpha\to0^+} \frac{f(t)}{\Gamma(\alpha+1)} = \frac{f(t)}{\Gamma(1)}=f(t),
$
while
\[
    |R_\alpha(t)|\leq \frac{\sup_{s\in[0,t]}|f'(s)|}{\Gamma(\alpha)} \int_0^t (t-\sigma)^\alpha d \sigma = \frac{\sup_{s\in[0,t]}|f'(s)|}{\Gamma(\alpha)(\alpha+1)}\stackrel{\alpha\to0^+}{\longrightarrow}0,
\]
where we used that $\Gamma(0^+)=\infty$ (see Figure \ref{figura:gamma}).
This proves identity \eqref{def:I0}.
It is worth mentioning that this result is true for less regular functions, but we will not show the proof to avoid entering more technical issues.
The composition formula $I_{\alpha_1} I_{\alpha_2} = I_{\alpha_1 + \alpha_2}$ can be derived using the properties of the Gamma function, but we will prove it in a simpler way below using another characterization of the fractional integral \eqref{eq:defintfrac}.

\paragraph{A ``modern'' presentation.}
Given a sufficiently regular function $g$, the Laplace transform $\Laplace$ of $g$ is defined as
\begin{equation}
    \Laplace[g](s)=\int _{0}^{\infty }g(t)\, e^{-st}\,dt.
\end{equation}
This operator is widely used due to its versatility for solving ordinary differential equations (both classical and fractional). In particular, the key property is
\begin{align*}
    \Laplace[I_1 f](s) & =  \int _{0}^{\infty }e^{-st} \int_0^t f(\sigma )\, d\sigma dt = \int _{0}^{\infty }f(\sigma ) \int_{\sigma}^\infty  e^{-st}\,dt  d\sigma \\
                       & =s^{-1} \int _{0}^{\infty }f(\sigma)\,  e^{-s\sigma} d\sigma  =s^{-1} \Laplace[f] (s).
\end{align*}
Because of this property, one says that $I_1$ has Laplace symbol $s^{-1}$. Iterating this trick we observe that, for $n\in \mathbb{N}$,
\[
    \Laplace[I_n f](s)=s^{-1}\Laplace[I_{n-1} f](s)= \cdots = s^{-n}\Laplace[f](s).
\]
This tells us that $I_n$ has Laplace symbol $s^{-n}$. Therefore, it is natural to define the fractional integral $I_\alpha$ with $\alpha\geq0$ as the operator with Laplace symbol $s^{-\alpha}$, that is, require that
\begin{equation}\label{eq:defintlap}
    \Laplace[I_\alpha f](s)=s^{-\alpha} \Laplace[f](s).
\end{equation}
In fact, some authors use this characterization as a definition. It is not difficult to check that the operator defined in \eqref{eq:defintfrac} satisfies \eqref{eq:defintlap}, and therefore both definitions coincide.

Moreover, as we announced, it is immediate to check the composition property
\begin{equation*}
    I_{\alpha_1} I_{\alpha_2} f=I_{\alpha_1+\alpha_2}f,
\end{equation*}
using \eqref{eq:defintlap}, since
\begin{align*}
    \Laplace[I_{\alpha_1}I_{\alpha_2} f](s) & =s^{-\alpha_1} \Laplace[I_{\alpha_2} f](s)=s^{-\alpha_1} s^{-\alpha_2} \Laplace[f](s)=s^{-\alpha_1-\alpha_2} \Laplace[f](s) \\
                                            & =\Laplace[I_{\alpha_1+\alpha_2} f](s).
\end{align*}

\subsection{Riemann–Liouville and Caputo fractional derivatives}\label{sec:defCapRL}

Once the way to integrate a non-integer number of times is understood, the idea of how to differentiate $\alpha\in \mathbb R_{>0} \setminus \mathbb N$ times arises naturally: just integrate $n-\alpha$ times and differentiate $n$ times (or vice versa) where $n=\lceil \alpha \rceil$.\footnote{$\lceil \alpha \rceil$ denotes the smallest natural number greater than or equal to $\alpha$} This gives rise to two (not equivalent!) definitions of fractional derivative.
Analogously to what we have said for integrals, the derivatives we introduce are well defined for functions $C^{\lceil \alpha \rceil} ([0,\infty))$, although the class of functions to which they can be applied is slightly broader. These details lie outside the scope of these notes.

The \textit{Riemann–Liouville derivative} is defined as
\begin{align}\label{eq:defC}
    \dRL f(t)\coloneqq \frac{d^{n}}{d t^{n}}[I_{n-\alpha} f](t)
    = \frac{1}{\Gamma(n-\alpha)}\frac{d^{n}}{d t^{n}} \int_{0}^t \frac{f(s)}{(t-s)^{\alpha+1-n}}ds,
\end{align}
while the \textit{Caputo derivative} is given by
\[
    \dC f(t)\coloneqq I_{n-\alpha} [f^{(n)}](t)= \frac{1}{\Gamma(n-\alpha)} \int_{0}^t \frac{f^{(n)}(s)}{(t-s)^{\alpha+1-n}}ds,
\]
where $f^{(n)}$ denotes the $n$-th derivative of $f$.
It is important to note at this point that fractional derivatives have the \emph{memory property} in the sense that they look at the entire ``previous history'' of $f$ (i.e. its values in $(0,t)$) to compute its fractional derivative at time $t$.

As we did with the integral, for $\alpha = 0$ we define
\[
    \dC[0] \coloneqq I_0\frac{d^0}{dx^0} = I_0, \qquad\dRL[0]  \coloneqq \frac{d^0}{dx^0} I_0 = I_0,
\]
where, as observed in  \eqref{def:I0}, $I_0$ is the identity operator $f \mapsto f$.

Recall that our initial goal was to define fractional derivatives so that they ``extend'' classical derivatives to non-integer indices. Let us see that, indeed, if $\alpha\in \mathbb{N}$ (and therefore $n=\lceil\alpha\rceil=\alpha$), then $\dC$ and $\dRL$ coincide with the usual derivative.
For the Riemann–Liouville derivative,
\begin{align*}
    \dRL f(t) & = \frac{d^n}{dt^n}I_{n-\alpha}[f](t) = \frac{d^n}{dt^n}I_{0}[f](t) = f^{(n)}(t)
    = f^{(\alpha)}(t).
\end{align*}
For the Caputo derivative,
\begin{align*}
    \dC f(t) & = I_{n-\alpha}[f^{(n)}](t) = I_{0}[f^{(n)}](t) = f^{(n)} (t) = f^{(\alpha)}(t).
\end{align*}

An important property (or rather limitation) is that, while ordinary derivatives are well defined at almost every point for any monomial $t^\beta$ with $\beta\in \mathbb{R}$, fractional derivatives are only defined up to a certain order. More precisely:
\begin{itemize}
    \item $\dC t^\beta$ is only defined if $\beta \ge 0$.
    \item $\dRL t^\beta$ is only defined if $\beta > -1$.
\end{itemize}
These ranges arise from the (non-)integrability at zero of the quantities involved in the definition.

As we have already mentioned, Caputo and Riemann–Liouville derivatives are not equivalent in general. In fact, the following relation holds for all $\alpha\in \mathbb{R}_{> 0}\setminus \mathbb{N}$:
\begin{align}
    \label{eq:relRL-C}
    \dRL f(t) = \dC f(t) + \sum_{k=0}^{\lceil \alpha \rceil - 1} \frac{f^{(k)}(0)}{\Gamma(k+1-\alpha)} t^{k-\alpha}.
\end{align}
We point out that the formula includes Gamma functions at negative values (see Section \ref{sec:gamma}).
Let us check \eqref{eq:relRL-C} for $\alpha \in (0,1)$ (and therefore $n=\lceil\alpha \rceil=1$) for simplicity. Integrating by parts and using Leibniz’s rule for differentiation under the integral we obtain
\begin{align*}
    \dRL f(t) & = \frac{1}{\Gamma(2-\alpha)} \frac{d}{dt} \left[ f(0) t^{1-\alpha} + \int_0^t f'(s)(t-s)^{1-\alpha} \, ds \right]                \\
              & = \frac{1-\alpha}{\Gamma(2-\alpha)}f(0)t^{-\alpha} + \frac{1-\alpha}{\Gamma(2-\alpha)}\int_0^t \frac{f'(s)}{(t-s)^{\alpha}}\, ds \\
              & =\frac{f(0)}{\Gamma(1-\alpha)}t^{-\alpha} +\dC f(t).
\end{align*}
Relation \eqref{eq:relRL-C} provides several interesting properties for $\alpha \in \mathbb R_{> 0} \setminus \mathbb  N$:
\begin{enumerate}[label=(\alph{enumi})]
    \item  The Riemann–Liouville and Caputo derivatives of a function $f$ coincide if and only if its first $\lceil \alpha \rceil - 1$ integer derivatives (understanding $f^{(0)}=f$) are zero at $t=0$.
    \item It is trivial to observe that for any constant function $g(t) = K$ we have
          \[
              \dC g(t) = 0,
          \]
          and therefore,
          \begin{equation}
              \label{eq:constante}
              \dRL g(t) = \frac{K}{\Gamma(1-\alpha)}t^{-\alpha}.
          \end{equation}
          In some sense, this tells us that the Riemann–Liouville derivative interprets the function $g$ as ``$g(t)=K t^{0}$'' and proceeds to reduce the order of the polynomial. From the definition \eqref{eq:defC} one can check (using \eqref{eq:Ialpha de monomio}) that, in fact, $\dRL$ reduces the order for general monomials of order $\beta>-1$:
          \begin{equation}
              \label{eq:der-monomio}
              \dRL t^{\beta}= \frac{\Gamma(\beta+1)}{\Gamma(\beta - \alpha +1)}  t^{\beta - \alpha}.
          \end{equation}
          Note also that $\dRL t^{\alpha - 1} = 0$.

    \item Another interesting property is that, for any function $f$ well defined on an interval $[0,T]$ and sufficiently regular (for example $|\frac{df}{dt} (t)| \le Ct^{\beta-1}$ for some $\beta > \alpha $), we have
          \[
              \lim_{t\to0^+} \dC f(t) = 0,
          \]
          while, if $f(0)\not=0$, then necessarily
          \[
              \lim_{t\to0^+} \left|\dRL f(t)\right| = \infty.
          \]
\end{enumerate}
\paragraph{Laplace transform of fractional derivatives.}
Recall that the Laplace transform of the derivative can also be computed easily by formally integrating by parts
\begin{align*}
    \Laplace\left[\frac{df}{dt} \right] (s) & = \int_0^\infty \frac{df}{dt} (t) e^{-st} dt
    =- \int_0^\infty f(t) \frac{d}{dt} (e^{-st}) dt - f(0^+)                               \\
                                            & = s \Laplace[f] (s) - f(0^+).
\end{align*}
Thus, for the Caputo derivative, using the Laplace transform properties of $\frac{d}{dt}$ and $I_\alpha$, if $\alpha \in (0,1)$ we have
\begin{equation}
    \label{eq:caputo-laplace}
    \begin{split}
        \Laplace \left[ \dC  f \right] (s) & = \Laplace \left[ I_{1-\alpha} \frac{d}{dt} f \right] (s)
        = s^{-(1-\alpha)} \Laplace\left[\frac{d}{dt} f\right]                                          \\
                                           & = s^\alpha \Laplace [f] - s^{\alpha - 1} f(0^+).
    \end{split}
\end{equation}
Similarly, for the Riemann–Liouville derivative, if $\alpha\in(0,1)$, then
\begin{equation*}
    \Laplace\left[ \dRL  f \right] (s) = s^\alpha \Laplace[f] (s) -   I_{1-\alpha} f(0^+).
\end{equation*}
Similar formulas can be written for $\alpha \in (1,2)$, $\alpha \in (2,3)$, etc., following the same procedure.
\normalcolor

Using this tool we can study when the semigroup property \eqref{eq:comp} holds for more general functions than monomials. We do it only in a simple case for brevity.
If $\alpha_1 , \alpha_2 > 0$ and $\alpha_1 + \alpha_2 < 1$, then
\begin{align*}
    \Laplace[\dC[\alpha_1] \dC[\alpha_2] f](s) & = s^{\alpha_1} \Laplace [\dC[\alpha_2] f](s) - s^{\alpha_1 - 1} \dC[\alpha_2] f(0^+)                                \\
                                               & = s^{\alpha_1 + \alpha_2} \Laplace [f](s) - s^{\alpha_1 + \alpha_2-1}f(0^+) - s^{\alpha_1 - 1} \dC[\alpha_2] f(0^+) \\
                                               & = \Laplace[\dC[\alpha_1 + \alpha_2] f](s) - s^{\alpha_1 - 1} \dC[\alpha_2] f(0^+).
\end{align*}
Thus, the semigroup property holds if and only if $\dC[\alpha_2] f (0^+) = 0$ (for which it suffices that $\left|\frac{df}{dt}\right| \le Ct^{\beta -1} $ with $\beta > \alpha_2$).
A similar statement holds for $\dRL$.

\paragraph{Relations with the usual derivative when $\alpha$ tends to a natural number.}
We have already seen that if $\alpha\in \mathbb{N}$ then $\dC$ and $\dRL$ coincide with the usual derivative.

It is also interesting to consider the continuity of the operators $\dC$ and $\dRL$ with respect to the parameter $\alpha$ when it tends to a natural number.
The integral operator $I_\alpha$ satisfies $I_\alpha f \to f$ as $\alpha \to 0^+$ and $I_\alpha f \to I f$ as $\alpha \to 1^-$ (the first was seen in subsection \ref{sec:intfrac} and the second is a direct consequence of the dominated convergence theorem). It is then not difficult to justify that, for sufficiently regular $f$, the following properties hold:
\begin{equation}\label{eq:limRL}
    \dRL f(t) = \frac{d}{dt} I_{1-\alpha} f (t)\longrightarrow
    \begin{dcases}
        f(t)              & \text{when } \alpha \to 0^+, \\
        \frac{df }{dt}(t) & \text{when } \alpha \to 1^-,
    \end{dcases}
\end{equation}
i.e., the Riemann–Liouville derivative tends to the identity when $\alpha\to0^+$ and to the usual derivative when $\alpha\to1^-$.
On the other hand, using the relation \eqref{eq:relRL-C} and the Gamma function properties $\Gamma(1) = 1$ and $\Gamma(0^+) = \infty$ (see \eqref{eq:gamma(0)}), we have
\begin{equation}\label{eq:limC}
    \dC f (t) = \dRL f (t) - \frac{f(0)}{\Gamma(1-\alpha)} t^{-\alpha}
    \longrightarrow
    \begin{dcases}
        f(t) - f(0)      & \text{when } \alpha \to 0^+, \\
        \frac{df}{dt}(t) & \text{when } \alpha \to 1^-.
    \end{dcases}
\end{equation}
To derive similar properties for higher-order derivatives it is helpful to note the following observation, which follows directly from the definitions of fractional derivatives:
\begin{equation*}
    \dRL[\alpha] f = \frac{d^m}{dt^m} \dRL[\alpha - m] f \quad  \text{ and } \quad   \dC[\alpha] f = \dC[\alpha - m] \frac{d^m}{dt^m} f, \qquad \forall m \in \mathbb N,\, m < \alpha.
\end{equation*}
Thus, from the limits as $\alpha \to 0^+, 1^-$ seen in \eqref{eq:limRL} and \eqref{eq:limC}, we deduce that for all $m \in \mathbb N$
\begin{equation*}
    \begin{aligned}
        \dRL f(t) & \to f^{(m)}(t) \quad \text{if } \alpha \to m,             \\
        \dC f(t)  & \to \begin{dcases}
                            f^{(m)} (t) - f^{(m)}(0) & \text{if } \alpha \to m^+, \\
                            f^{(m)} (t)              & \text{if } \alpha \to m^-.
                        \end{dcases}
    \end{aligned}
\end{equation*}
Note that the Caputo derivative can be discontinuous with respect to $\alpha$ at the values $\alpha \in \mathbb N$, and the left limit coincides with the pointwise value $\dC[m] f (t) = f^{(m)}(t)$.
This equals the right limit if and only if $f^{(m)} (0) = 0$.

\subsection{The Grünwald-Letnikov fractional derivative}

This formalization, based on the finite-difference definition of higher-order derivatives, was introduced by A. K. Grünwald in Prague in 1867 \cite{Grunwald1867}, and by A. V. Letnikov in Moscow in 1868 \cite{Letnikov1868}.
We start from the definition of derivative by finite differences
\[
    f^{(1)}(t)= \lim_{h\to0^+} \frac{f(t)-f(t-h)}{h}.
\]
Higher-order derivatives can be obtained by iterating to get
\begin{align*}
    f^{(2)}(t) & = \lim_{h\to0^+} \frac{f^{(1)}(t)-f^{(1)}(t-h)}{h}         \\
               & =\lim_{h\to0^+} \frac{(f(t)-f(t-h))-(f(t-h)-f(t-2h))}{h^2} \\
               & = \lim_{h\to0^+} \frac{f(t)-2f(t-h)+f(t-2h)}{h^2},
\end{align*}
and, more generally, for any $n\in \mathbb{N}$, the formula
\begin{align*}
    f^{(n)}(t) & = \lim_{h\to0^+} \frac{1}{h^n} \sum_{k=0}^n(-1)^k \binom{n}{k} f(t-kh).
\end{align*}
We note that
\begin{equation*}
    \binom{n}{k} =\frac{\Gamma(n+1)}{\Gamma(k+1) \Gamma(n-k+1)}.
\end{equation*}
In this way, it is natural to define the \textit{Grünwald-Letnikov derivative} of order $\alpha\geq0$ as follows:
\begin{align*}
    \dGL(t) & = \lim_{j\to+\infty} \frac{1}{(h_j)^\alpha} \sum_{k=0}^{j}(-1)^k \frac{\Gamma(\alpha+1)}{\Gamma(k+1) \Gamma(\alpha-k+1)}f(t-kh_j), \quad \text{where } j h_j=t.
\end{align*}

We now make some interesting observations about this concept of derivative:
\begin{enumerate}[label=(\alph{enumi}), align=left, leftmargin=*]
    \item
          The Grünwald-Letnikov derivative $\dGL$ coincides with the classical derivative $\frac{d^{(\alpha)}}{dt^{(\alpha)}}$ if $\alpha=n\in \mathbb{N}$, since $|\Gamma(-m)|= \infty$ when $m\in\mathbb{N}$ (see Section \ref{sec:gamma}), which cancels the terms with $k \geq n+1$.

    \item
          It is an interesting (and nontrivial) exercise to show that
          the Riemann-Liouville and Grünwald-Letnikov derivatives coincide, that is,
          \[
              \dGL = \dRL.
          \]
          The proof is based on interpreting the formula for $\dGL$ as a numerical quadrature of the formula for the Caputo derivative given in \eqref{eq:defC} plus a remainder, and then using the relation between Caputo and Riemann-Liouville derivatives given in \eqref{eq:relRL-C}.

    \item It is important to keep in mind that the memory property is observed in the Grünwald-Letnikov derivative in the fact that, again, to compute the derivative at time $t$ one looks at the entire previous history of the function in $[0,t]$. This is because, even before taking the limit, $f$ is evaluated on finer and finer equally spaced meshes of $[0,t]$.

    \item In view of the previous comments, the formulation of $\dGL$ naturally suggests a numerical discretization of the Riemann-Liouville and Caputo derivatives. However, in practice, other methods are used that offer better results (see Section \ref{sec:metodos numericos}).
\end{enumerate}

\section{Fundamental theorems of fractional calculus} \label{sec: teoremas fundamentales}

The fundamental theorem of (integer) calculus can be interpreted in two ways. The first, and most usual, is that
\[
    f(t) = \frac{d}{dt} \int_0^t f(s) ds = \partial_t I f (t).
\]
The second interpretation, also called Barrow’s rule, tells us that
\[
    f(t) - f(0) = \int_0^t f' (s) ds = I \partial_t f (t) .
\]
These results, which have innumerable applications in the daily life of any mathematician, have their corresponding fractional versions that we present below.

First let us try to reproduce the fundamental theorem of calculus in the fractional case with $\alpha\in(0,1)$.
For $\dRL$ we can use the semigroup property of $I_\alpha$ and deduce
\begin{equation}
    \label{eq:dRL Ialpha}
    \dRL I_{\alpha} f(t) =\partial_t I_{1-\alpha} I_\alpha f(t) = \partial_t I f (t) = f(t).
\end{equation}
For $\dC$, we use its relation with $\dRL$ given in \eqref{eq:relRL-C}, to observe that
\begin{equation}
    \label{eq:dC Ialpha}
    \begin{split}
        \dC I_\alpha  f (t) & = \dRL I_\alpha  f (t) - \frac{t^{-\alpha}}{\Gamma(1-\alpha)} I_{\alpha} f (0^+) \\
                            & = f(t) - \frac{t^{-\alpha}}{\Gamma(1-\alpha)} I_{\alpha} f (0^+) .
    \end{split}
\end{equation}
Note that the “additional term” $\frac{t^{-\alpha}}{\Gamma(1-\alpha)} I_{\alpha} f (0^+)$ vanishes, for example, if $f$ is bounded in a neighborhood of $0$.

Let us now proceed with Barrow’s rule.
On the one hand we have, again by the semigroup property of the fractional integral, that
\begin{equation}
    \label{eq:Ialpha dC}
    I_{\alpha} \dC f(t) = I_\alpha I_{1-\alpha} \partial_t f (t) = I \partial_t f (t) = f(t)-f(0).
\end{equation}
On the other hand, using the fundamental theorem of fractional calculus \eqref{eq:dC Ialpha}, we have that
\begin{equation}
    \begin{split}
        \label{eq:Ialpha dRL}
        I_\alpha \dRL f (t) & =I_\alpha \frac{d}{dt} I_{1-\alpha}f(t)= \dC[1-\alpha]I_{1-\alpha}f(t) \\
                            & = f(t) - \frac{t^{\alpha-1}}{\Gamma(\alpha)} I_{1-\alpha} f (0^+).
    \end{split}
\end{equation}
Again, if $f$ is bounded at $0$ the “additional” term vanishes.
Note that, when $\alpha \to 1^-$, the usual formulas are recovered.

We can rewrite the Barrow rules to recover a function from its derivatives and “information” at $t = 0$. In particular, if $\dC u = g_1$ and $\dRL v = g_2$ for certain given functions $g_1$ and $g_2$, then we can recover $u$ and $v$ using \eqref{eq:Ialpha dC} and \eqref{eq:Ialpha dRL}. More precisely,
\begin{align*}
     & u(t)-u(0)=I_\alpha \dC u(t)= I_\alpha g_1(t),                                                          \\
     & v(t)- \frac{t^{\alpha-1}}{\Gamma(\alpha)} I_{1-\alpha} v (0^+) = I_\alpha \dRL v(t) = I_\alpha g_2(t),
\end{align*}
which we can rewrite, using the usual notation, as
\begin{align}
    \label{eq:u from dC u}
    u(t) & = u(0) + \int_0^t g_1(s) K_{\alpha} (t-s) ds                             \\
    \label{eq:v from dRL v}
    v(t) & = I_{1-\alpha} v(0^+) K_\alpha(t) + \int_0^t g_2(s) K_{\alpha} (t-s) ds.
\end{align}
These formulations will be useful in the next section, since they will allow us to express initial value problems in integral form.

\section{Ordinary Differential Equations}\label{sec:ODE}

The aim of this section is to present some basic techniques to solve linear and nonlinear ordinary differential equations. We start by reviewing classical techniques for local problems, and then see how to extend the reasoning to nonlocal problems.

\subsection{Local equations: linear homogeneous case}
\label{sec:edos locales}
The linear ordinary differential equation
\begin{equation}
    \label{eq:EDO local lineal homogenea}
    \frac{du}{dt} = \lambda u
\end{equation}
models the evolution of a quantity $u$ when the variation of $u$ (its derivative) is proportional to the quantity $u$.
This model is known as the \emph{Malthusian model}.
It is a simple model used in different fields: from radioactive decay (where $\lambda < 0$ is related to the half-life) to population dynamics (where $\lambda$ is related to the birth/death rate).

A classical way to look for solutions to this equation is to study the power series of a possible regular solution
\[
    u(t) = a_0 + a_1 t + a_2 t^2 + \cdots + a_n t^n + \cdots
\]
To find the coefficients we can formally substitute $u$ into the differential equation \eqref{eq:EDO local lineal homogenea}, obtaining
\[
    a_1 + 2 a_2 t
    + \cdots
    + n a_n t^{n-1}
    + \cdots = \lambda (a_0 + a_1 t + a_2 t^2 + \cdots + a_n t^n + \cdots ).
\]
For the above identity to hold we must have $n a_{n} = \lambda a_{n-1}$ for every $n \ge 1$.
We can use finite recursion to see that $a_n = a_0 \lambda^n / n!$.
These are the coefficients of the well-known series of the exponential function, that is,
\[
    u(t)
    = a_0\sum_{n=0}^\infty\frac{(\lambda t)^n}{n!}
    = a_0 e^{\lambda t}.
\]
Moreover, we see that $u(0)=a_0$. This function $u$ is the unique solution of
\[
    \begin{dcases}
        \frac{du}{dt} (t) = \lambda u (t) & \text{for all } t > 0, \\
        u(0) = a_0.
    \end{dcases}
\]
This is known as an \textit{initial value problem}: given the initial state $u(0)$ and the model (the differential equation) we can know the state of the system at any later time.

\subsection{Fractional equations: linear homogeneous case}
Let us now see what happens when we replace the classical derivative by the Caputo and Riemann–Liouville derivatives. For simplicity, we consider only the case $\alpha\in (0,1)$.

\paragraph{With Caputo derivative.}
We first study the equation
\[
    \dC u = \lambda u.
\]

Functions that coincide with their power series are, in general, very regular. This is the case of the exponential function. Solutions of fractional equations are not so regular and require a bit more analysis.
We provide a heuristic argument for this fact.
Looking at Barrow’s rule \eqref{eq:u from dC u} with $g_1(t) = \lambda u(t)$, if we assume $u(0) > 0$ and $u$ is continuous in a neighborhood of zero then $u(t) = u(0) + \lambda I_\alpha u (t)$.
For $t$ sufficiently small such that $u(s) > 0$ for $s \in [0,t]$ (we can assume the existence of such $t$ due to continuity) we have
\begin{align*}
    I_\alpha u (t) & \ge \min_{s\in[0,t]} u(s) \int_0^t K_\alpha (t-s) ds = \frac{t^\alpha}{\Gamma(1+\alpha)} \min_{s\in[0,t]} u(s), \\
    I_\alpha u (t) & \le \max_{s\in[0,t]} u(s) \int_0^t K_\alpha (t-s) ds = \frac{t^\alpha}{\Gamma(1+\alpha)} \max_{s\in[0,t]} u(s).
\end{align*}
We deduce that $u\sim t^{\alpha}$ for small $t$,\footnote{We will write $f\sim g$ to indicate that the functions $f$ and $g$ “behave approximately the same”, without entering into details.}
and therefore $u$ is not very regular (it is not even $C^1$), so it is not reasonable to study the power series of $u$.
However, following the previous heuristic reasoning, it seems reasonable to look for $u(t) = U(t^\alpha)$ with $U$ regular. Using the power series of $U$, we try to write $u$ as
\[
    u(t) = a_0 + a_1 t^\alpha + a_2 t^{2\alpha} + \cdots + a_n t^{n\alpha} + \cdots
\]
For $\alpha = 1$ we saw by recursion that $a_n = a_0 \lambda^n/n! = a_0 \lambda^n/\Gamma(n+1)$. With similar (though more tedious) arguments, one proves that in this case the solution is $a_n = a_0 \lambda^n / \Gamma(\alpha n + 1)$. Therefore, the solution can be written as
\begin{equation}\label{eq:ML}
    u(t) = a_0 E_{\alpha,1} (\lambda t^\alpha), \qquad \text{where } E_{\alpha,\beta} (z) = \sum_{n=0}^\infty \frac{z^n}{\Gamma(\alpha n+\beta)}.
\end{equation}
The functions $E_{\alpha, \beta}$ are known as \textit{Mittag-Leffler functions}, and
it is common to denote $E_\alpha (z) \coloneqq E_{\alpha,1} (z)$.
This function is continuous and satisfies $E_{\alpha} (0) = 1$, so again we have $u(0) = a_0$. Thus, $u$ is the unique solution of the initial value problem
\[
    \begin{dcases}
        \dC u (t) = \lambda u (t) & \text{for all } t > 0, \\
        u(0) = a_0.
    \end{dcases}
\]
We note that the solutions are not regular beyond $C^\alpha$. In fact, one can show that $u^{(m)}(t)\sim t^{\alpha-m}$ as $t\to 0^+$ for all $m\in \mathbb{N}$.

\paragraph{With Riemann–Liouville derivative.}
We now seek solutions of
\begin{equation}
    \label{eq:EDO RL lineal homogeneo}
    \dRL v = \lambda v.
\end{equation}
As we saw in \eqref{eq:constante},
$\dRL 1 = \frac{t^{-\alpha}}{\Gamma(1-\alpha)}$, but the power series expansion of regular functions has no negative exponent terms.
We cannot expect to find solutions $v(t) = V(t^\alpha)$ either.

It is not easy to guess what the structure of the solution of \eqref{eq:EDO RL lineal homogeneo} will be (we will see that, in fact, it is $t^{\alpha - 1} V(t^\alpha)$).
However, we can relate this problem for Riemann–Liouville derivatives with the one presented above for the Caputo derivative.
Applying $I_{1-\alpha}$ to equation \eqref{eq:EDO RL lineal homogeneo} we arrive at
\[
    I_{1-\alpha}\frac{d}{dt} I_{1-\alpha} v = \lambda I_{1-\alpha} v.
\]
Setting $u = I_{1-\alpha} v$ we have $\dC u = \lambda u$, hence $u(t) = a_0 E_{\alpha}(\lambda t^\alpha)$, where $a_0 = u(0^+)=I_{1-\alpha} v (0^+)$.
Now, using the fundamental theorem of fractional calculus \eqref{eq:dRL Ialpha} and formula \eqref{eq:der-monomio}, we obtain
\begin{align*}
    v(t) & = \dRL[1-\alpha] I_{1-\alpha} v (t) = \dRL[1-\alpha] u(t)                                                                                                 \\
         & = a_0 \sum_{n=0}^\infty \frac{\dRL[1-\alpha](\lambda^n  t^{\alpha n}) }{\Gamma(\alpha n + 1)}                                                             \\
         & = a_0 \sum_{n=0}^\infty \frac{\lambda^n  \frac{\Gamma(\alpha n + 1)}{\Gamma(\alpha n + 1 - (1-\alpha))} t^{\alpha n - (1-\alpha)} }{\Gamma(\alpha n + 1)} \\
         & = a_0 t^{\alpha - 1} E_{\alpha,\alpha} (\lambda t^\alpha).
\end{align*}
We see that, if $a_0 \ne 0$, then $|v(0^+)| = \infty$. This is surprising.

Let us analyze what kind of “information” at time $t = 0$ is natural for this type of problems.
From how we defined $u$ in terms of $v$, we observe that
\[
    a_0 = u(0^+) = I_{1-\alpha} v (0^+).
\]
As in Barrow’s rule \eqref{eq:v from dRL v}, we see that the natural “initial condition” for this type of problems is not the value $v(0^+)$, but the singular initial “datum” $I_{1-\alpha} v (0^+)$.
Moreover, $v$ is the unique solution of the singular initial value problem
\[
    \begin{dcases}
        \dRL v(t) = \lambda v(t) & \text{for all } t > 0, \\
        I_{1-\alpha} v(0^+) = a_0.
    \end{dcases}
\]
Note that this problem does not admit finite $v(0^+)$ conditions different from $0$.
Indeed, if $|v(0^+)| < \infty$, then $a_0 = 0$, and $v(t) = 0$ for all $t > 0$.

It is customary to denote
\begin{equation}
    \label{eq:Palpha}
    P_{\alpha}(t;-\lambda) \coloneqq t^{\alpha-1} E_{\alpha,\alpha}(\lambda t^\alpha).
\end{equation}
This notation is due to the fact that, in some contexts, it is natural that $\lambda < 0$.

\subsection{Non-homogeneous linear problems}
\label{sec:EDO lineal no homogenea por Laplace}
We now turn to non-homogeneous linear problems, that is, those in which the right-hand side of the differential equation depends not only on the solution but is also influenced by external “forces” that may vary with time.
We study only the case in which the coefficients of the linear part are constant.
Although in some cases we could continue arguing with power series (as in the previous sections), their generalization to the fractional setting is better treated through the Laplace transform, so we proceed in this way.

\paragraph{Local problems.}
We consider the initial value problem
\[
    \begin{dcases}
        \frac{du}{dt}(t) = \lambda u(t) + f(t) & \text{for all } t > 0, \\
        u(0) = u_0,
    \end{dcases}
\]
where $f$ is a given function.
Applying the Laplace transform to both sides of the differential equation we obtain
\[
    s \Laplace[u] (s) - u_0 = \lambda \Laplace[u](s) + \Laplace[f] (s),
\]
which, after solving for $\Laplace[u](s)$, yields a closed formula for the transform of the solution $u$ given by
\[
    \Laplace[u](s) =  \frac{u_0}{s-\lambda}   + \frac{\Laplace[f](s)}{s-\lambda}.
\]
On the one hand, it is immediate to check that
for  $g(t)=e^{\lambda t}$ we have $\Laplace[g] (s) = \frac{1}{s-\lambda}$.
On the other hand, using the Laplace \emph{convolution} property
\begin{equation*}
    \label{eq:Laplace convolution}
    \begin{aligned}
        \Laplace\left[ \int_0^\cdot f(\sigma) g(\cdot-\sigma) d\sigma \right] (s)
         & = \Laplace[f] (s) \Laplace[g](s),
    \end{aligned}
\end{equation*}
we deduce that the solution of our problem is given by
\[
    u(t) = u_0 e^{\lambda t} + \int_0^t e^{\lambda(t-s)} f(s) ds.
\]
This is the so-called \emph{variation of constants formula} in the local case.

\paragraph{Fractional problems.} We now consider the problem
\begin{equation}
    \label{eq:ODE Caputo}
    \begin{dcases}
        \dC u (t) = \lambda u (t) + f (t) & \text{for all } t > 0 \\
        u(0) = u_0.
    \end{dcases}
\end{equation}
Following the same reasoning as in the local case, and using property \eqref{eq:caputo-laplace} of the Laplace transform applied to the Caputo derivative, we obtain
\[
    \Laplace [u] (s)
    = \frac{s^{\alpha-1}}{s^\alpha - \lambda} u_0 + \frac{1}{s^\alpha - \lambda} \Laplace[f](s).
\]
Fortunately, the Mittag-Leffler functions have the following Laplace transform property:
\begin{equation}
    \label{eq:Mittag-Leffler Laplace}
    \begin{aligned}
        \Laplace[t^{\beta-1}E_{\alpha,\beta}(\lambda t^\alpha)](s)
         & =
        \int_0^\infty
        t^{\beta - 1} E_{\alpha, \beta} (\lambda t^\alpha) e^{-st} d t                   \\
         & = \sum_{k=0}^\infty \frac{\lambda^k}{\Gamma(\alpha k + \beta)}  \int_0^\infty
        t^{\beta + \alpha k - 1}  e^{-st} d t                                            \\
         & = \sum_{k=0}^\infty \frac{\lambda^k}{s^{k\alpha + \beta}}
        = s^{-\beta}\sum_{k=0}^\infty \left( \frac{\lambda}{s^\alpha} \right)^k
        = s^{-\beta} \frac{1}{1-\frac \lambda {s^\alpha}}                                \\
         & =\frac{s^{\alpha - \beta}}{s^\alpha - \lambda},
    \end{aligned}
\end{equation}
when $\alpha, \beta > 0$.
Thus we rewrite the transform of $u$ as
$$
    \begin{aligned}
        \Laplace [u] (s)
         & =u_0 \Laplace[E_\alpha(\lambda t^\alpha)](s)
        +\Laplace[t^{\alpha-1}E_{\alpha,\alpha}(\lambda t^\alpha)](s)
        \Laplace[f](s).
    \end{aligned}
$$
Using the notation \eqref{eq:Palpha} we deduce that the general solution of \eqref{eq:ODE Caputo} is
\begin{equation}
    \label{eq:ODE Caputo solution}
    u (t) = u_0 E_\alpha (\lambda t^\alpha) + \int_0^t P_{\alpha}(t-\tau;-\lambda) f(\tau) d \tau.
\end{equation}

On the other hand, the equation
\begin{gather}
    \label{eq:ODE Riemann}
    \begin{dcases}
        \dRL v(t)  = \lambda v(t)  +  f(t) , & \text{for all }  t > 0 , \\
        I_{1-\alpha} v (0^+) = v_0,
    \end{dcases}
\end{gather}
becomes
\begin{equation*}
    \Laplace [v] (s) = \frac{1}{ s^\alpha - \lambda}v_0 + \frac{1}{s^\alpha - \lambda}\Laplace[f](s),
\end{equation*}
and therefore the general solution of \eqref{eq:ODE Riemann} is
\begin{equation}
    \label{eq:ODE Riemann solution}
    v(t) =   v_0 P_{\alpha}(t;-\lambda)+ \int_0^t P_{\alpha} (t-\tau;-\lambda) f(\tau) d \tau .
\end{equation}
We have obtained the variation of constants formulas in the fractional case, given by  \eqref{eq:ODE Caputo solution} and \eqref{eq:ODE Riemann solution}.

\subsection{Nonlinear problems}
\label{sec:no lineal. nocion de solucion}

In the local case, the nonlinear equation $u' = f(t,u)$ does not, in general, have such a simple closed-form solution. To prove properties of the solution one usually passes to the integral formulation (based on Barrow’s rule)
\begin{equation}
    \label{eq:edo local form integral}
    u(t) = u(0) + \int_0^t f(s, u(s)) ds, \qquad \text{for all } t \ge 0.
\end{equation}
When $f$ is regular, \eqref{eq:edo local form integral} has continuous solutions found by applying a Banach fixed point argument.

In the fractional case, we can apply the fractional Barrow rules \eqref{eq:u from dC u} and \eqref{eq:v from dRL v} to the case in which $g_1(t)=f_1(t,u(t))$ and $g_2(t)=f_2(t,u(t))$, and we obtain the integral formulations of the fractional ODEs
\begin{align}
    \label{eq:edo Caputo form integral}
    u(t) & = u(0) + \int_0^t f_1(s,u(s)) K_\alpha(t-s) ds ,                           \\
    \label{eq:edo RL form integral}
    v(t) & = K_\alpha(t) I_{1-\alpha} v(0^+) + \int_0^t f_2(s,v(s)) K_\alpha(t-s) ds.
\end{align}
As in the previous case, existence can be proved by techniques similar to the local case.

\section{Numerical methods}\label{sec:metodos numericos}

Except in very simple cases, the solution of a nonlinear problem usually cannot be expressed in terms of elementary functions. However, because of its usefulness in the real world, we would like to compute a suitable approximation that allows us to make decisions. This is where numerical methods appear as an indispensable tool. The goal of this section is to review some basic concepts of the local case and see how they extend to the fractional case.

\subsection{Local case}

From the very definition, several approximations of the derivative arise naturally. The most common ones, due to their numerical stability properties, are given by
\[
    u'(t) \approx \frac{u(t+\tau) - u(t)}{\tau}, \qquad u'(t) \approx  \frac{u(t) - u(t-\tau)}{\tau},
\]
for a chosen “time step” $\tau>0$.
This allows us to find the following approximate problems for differential equations of the form  $u'(t) = f(t,u(t))$:
\[
    \frac{u(t+\tau) - u(t)}{\tau} \approx u'(t) = f(t,u(t)), \quad \frac{u(t+\tau) - u(t)}{\tau} \approx u'(t+\tau) = f(t+\tau,u(t+\tau)).
\]
These are the so-called forward and backward approximations.
If we start at time $t = 0$ and define $t_n = n\tau$ with $n\in\mathbb{N}$, we can then think that $u(t_n) \approx U_n$, where $U_n$ solves the recurrence formulas (sometimes called difference equations)
\[
    \frac{U_{n+1} - U_n}{\tau} = f(t_n,U_n), \qquad \frac{U_{n+1} - U_{n}}{\tau} = f(t_{n+1}, U_{n+1}),
\]
with $U^0=u(0)$. If we construct $U_1, U_2, \ldots$ from these recurrence formulas, we arrive at the so-called explicit (forward Euler) and implicit (backward Euler) Euler methods. It is common to write them grouping present and future terms as follows:
\[
    U_{n+1} = U_n + \tau f(t_n,U_n), \qquad U_{n+1} - {\tau} f(t_{n+1}, U_{n+1}) = U_{n}.
\]

These same numerical methods can also be derived from the integral formulation \eqref{eq:edo local form integral} of the initial value problem (with the idea of applying a similar approach to the fractional case).
Suppose we know the values $U_0, \ldots, U_{n-1}$ and want to compute $U_{n}$.
The idea is to try to replace $f(t, u(t))$ with a function that uses only these data, for which there exist various reasonable approximations.
For example, the piecewise constant explicit and implicit interpolations
\[
    F_e(t)=f(t_j, U_j) \quad \textup{y} \quad F_i(t)=f(t_{j+1}, U_{j+1}) \quad \text{si} \quad t \in [t_j, t_{j+1}) \quad \textup{for} \quad j\in \mathbb{N}.
\]
The approximation $f(t,u(t)) \approx F_e(t)$ gives us
\begin{align*}
    U_{n+1} & \approx u(t_{n+1}) = u(0) + \int_0^t f(s,u(s)) ds \approx u(0) + \int_0^t F_e (s) ds \\
            & = u(0) + \tau \sum_{k=0}^n f(t_n,U_n).
\end{align*}
Thus, constructing by recurrence the sequence
\begin{equation*}
    U_{n+1} \coloneqq u(0) + \tau \sum_{k=0}^n f(t_n,U_n)
\end{equation*}
we observe that $U_{n+1} - U_n = \tau f(t_n,U_n)$, which is precisely the explicit Euler method.
Being able to make this simplification clearly shows that we are dealing with a local problem.
Similarly, using $F_i$ we obtain the implicit Euler method.
We could also use the linear interpolation
\[
    F_{\ell}(t) = f(t_k,u_k)\frac{t_{k+1} - t}{\tau} + f(t_{k+1},u_{k+1}) \frac{t-t_k}{\tau} \qquad \text{for } t \in [t_k, t_{k+1}],
\]
or even higher-order interpolations. This gives rise to the so-called Adams methods.

\subsection{Fractional case}

If we want to approximate the Caputo and Riemann-Liouville derivatives we can employ the discretization of the classical derivative. To do so, first decompose
\begin{align*}
    \dC u (t_n) & =  \int_0^{t_n} K_{1-\alpha} (t_n-s) \frac{d u}{dt} (s) ds                      \\
                & = \sum_{m=1}^n \int_{t_{m-1}}^{t_m} K_{1-\alpha} (t_n-s) \frac{d u}{dt} (s) ds.
\end{align*}
We can now leverage some of the discretizations of the usual derivative, for example\footnote{In this case we have used a piecewise linear interpolation of $u$ on each interval $[t_{m},t_{m+1}]$.}
\[
    \dC u (t_n) \approx \sum_{m=1}^n \frac{u(t_m) - u(t_{m-1})}{\tau} \int_{t_{m-1}}^{t_m} K_{1-\alpha} (t_n-s) {ds}.
\]
These latter integrals do not depend on $u$, and are reasonably easy to compute. From here we could obtain the explicit and implicit methods for an ordinary equation, but they are better derived from the integral formulation as we explain later.

For a regular function $u$, one could try to discretize $\dRL$ using this discretization of $\dC$ and \eqref{eq:relRL-C}.
This presents some subtleties, since the solution of \eqref{eq:ODE Riemann} has $|v(0^+)| = \infty$.
We will not go into further details here to avoid technical issues.

Recall that the fractional versions of the integral formulation \eqref{eq:edo local form integral} of a local initial value problem are given by \eqref{eq:edo Caputo form integral} and \eqref{eq:edo RL form integral}.
We can use ideas similar to the local case.
Fortunately, for both problems with $\dC$ and with $\dRL$ we need to find a discrete version of the following integral expression:
\[
    \int_0^t f(s,u(s)) K_\alpha(t-s)  ds.
\]
Thus, we can build Adams methods in the fractional case.
For $\dC$ we can use $ F_e,  F_i,  F_\ell$.
For $\dRL$, since $|v_0|=\infty$, at least at the first time step we can only use $ F_i$.
If $F(t)$ is constant on $[t_n,t_{n+1})$ for all $n\in \mathbb{N}$ (for example $ F=F_e$ or $F=F_i$), then
\begin{align*}
    \int_0^{t_{n+1}} F(t) K_\alpha(t_{n+1} - t) dt & = \sum_{k=0}^{n} F(t_k) \int_{t_k}^{t_{k+1}} K_\alpha(t_{n+1} - t) dt \eqqcolon  \sum_{k=0}^{n} F(t_k) w_{n+1,k}.
\end{align*}
The values of the integrals $w_{n,k}$ can, again, be computed explicitly, and do not depend on $F$. This produces the explicit and implicit Euler methods for the problem with Caputo derivative:
\begin{align*}
    U_{n+1} & = U_0 + \sum_{k=0}^{n} f(t_k, U_k) w_{n+1,k},          \\
    U_{n+1} & = U_0 + \sum_{k=0}^{n} f(t_{k+1}, U_{k+1}) w_{n+1,k} .
\end{align*}
Note that, unlike the local case, we cannot deduce $U_{n+1}$ solely from $U_n$.
In fact, to compute any step, we need all the previous ones. Here we see, again, the nonlocality of the fractional case.
In this direction different works have been developed; see for example \cite{DiethelmFordFreed2004DetailedErrorAnalysis} and its references.
\section{Systems of differential equations}\label{sec:sistemas}
To complete the study of the theory of fractional derivatives, and before presenting a couple of applications, we make some comments about what happens with systems of differential equations. As in previous sections, we start by reviewing the local case and we then see how similar ideas work in the fractional case.

\subsection{Local case}
Consider a system of equations of the form
\begin{equation*}
    \begin{dcases}
        \frac{d u_1}{dt} (t) = f_1(t,u_1(t), \cdots, u_n(t)) , \\
        \frac{d u_2}{dt} (t) = f_2(t,u_1(t), \cdots, u_n(t)) , \\
        \hspace{11ex} \vdots                                   \\
        \frac{d u_n}{dt} (t) = f_n(t,u_1(t), \cdots, u_n(t)) , \\
    \end{dcases}
\end{equation*}
which we write for convenience in its vector form
\begin{equation}
    \label{eq:sistema local no-lin}
    \frac{d \mathbf u}{dt} (t) = \mathbf f (t, \mathbf u(t)).
\end{equation}
Existence (and uniqueness) of solutions follows, as in the case of ordinary equations, by looking for a fixed point of the integral equation
\[
    \mathbf u(t) = \mathbf u_0 + \int_0^t \mathbf f(s, \mathbf u(s)) ds.
\]

\paragraph{Linear homogeneous problems.}

First assume that $\mathbf f(t,\mathbf u) = \mathbf A \mathbf u$ where $\mathbf A$ is a constant coefficient matrix (it does not change with time), that is,
\begin{equation}
    \label{eq:sistema local lineal}
    \frac{d \mathbf u}{dt}(t) = \mathbf A \mathbf u(t).
\end{equation}
As usual in the linear world, the solution $\mathbf u^\star(t) = \mathbf 0$ is distinguished.
One of the main points of interest is to know whether $\mathbf 0$ is a global attractor, that is, whether for any other solution $\mathbf u(t)$ with $\mathbf u(0) \ne \mathbf 0$ we have $\mathbf u(t) \to \mathbf 0$ as $t \to \infty$.

If $\mathbf w$ is an eigenvector of $\mathbf A$ with eigenvalue $\lambda \in \mathbb C$ (i.e. $\mathbf A\mathbf w= \lambda \mathbf w$), then one can find solutions of the form $\mathbf u(t) = c(t) \mathbf w$ for a certain function $c(t)$. Substituting into \eqref{eq:sistema local lineal} we have
\[
    c'(t) \mathbf w = c(t) \mathbf A\mathbf w= \lambda c(t) \mathbf w.
\]
This tells us that the function $c(t)$ must satisfy the linear ordinary differential equation $c' = \lambda c$ which, as we have already seen, has solution $c(t) = c(0) e^{\lambda t}$.
Note that, in this case, $\lambda$ may be a complex number. Therefore, it is pertinent to recall that the formula for the exponential of a complex number $\lambda = a + i b$ with $a,b \in \mathbb R$ is given by
\[
    e^{\lambda t} = e^{a t} (\cos(b t) + i\sin(b t)).
\]
This suggests that, if there exists an eigenvalue $\lambda$ with $a = \Re(\lambda) > 0$, then along the direction of $\mathbf{w}$ we tend to move away from $\mathbf 0$.
On the contrary, if $\Re(\lambda) < 0$ we approach $\mathbf 0$ along the line with direction vector $\mathbf{w}$.

Recalling the basic differential equations course (see, for example, \cite{Simmons}),
we know that the solution of problem \eqref{eq:sistema local lineal} for a general datum $\mathbf u(0) \in \mathbb R^n$ can be computed through the Jordan factorization of $\mathbf A$, and is
\[
    \mathbf u(t) = \sum_{j=1}^n e^{\lambda_j t} P_j(t) \mathbf w_j,
\]
where $P_j(t)$ are polynomials and $\mathbf w_j$ vectors that depend on $\mathbf A$ and $\mathbf u(0)$.
If the matrix $\mathbf A$ is diagonalizable, the polynomials reduce to constants.
Thus, a sufficient condition for $\mathbf 0$ to be an attractor of system \eqref{eq:sistema local linealizado} is that all eigenvalues $\lambda_j=a_j+i b_j$ of $\mathbf A$ (possibly complex numbers) satisfy $a_j = \Re (\lambda_j) < 0$.

\paragraph{Stability in autonomous nonlinear problems.}
Now assume that $\mathbf f(t,\mathbf u) = \mathbf f(\mathbf u)$, which is known as autonomous systems.
As in the linear case, steady states are constant solutions, that is, if $\mathbf u^\star \in \mathbb R^n$ is such that
\begin{equation*}
    \mathbf f(\mathbf u^\star) = \mathbf 0 ,
\end{equation*}
then $\mathbf u^\star (t) = \mathbf u^\star$ is a solution of the ODE.
To study stability of these constant solutions, we approximate $\mathbf f$ by its Taylor expansion
\[
    \mathbf f ( \mathbf u^\star + \mathbf z ) = \mathbf f (\mathbf u^\star) + \mathbf J\mathbf f ( \mathbf u^\star ) \mathbf z  + O(|\mathbf z|^2) = \mathbf J\mathbf f ( \mathbf u^\star ) \mathbf z + O(|\mathbf z|^2),
\]
where $\mathbf J \mathbf f$ denotes the Jacobian of $\mathbf f$.
If we decompose the solutions in the form $\mathbf u(t) = \mathbf u^\star + \mathbf z(t)$ we find that
\[
    \frac{d \mathbf z}{dt} (t) = \mathbf J\mathbf f ( \mathbf u^\star ) \mathbf z (t) + O(|\mathbf z (t)|^2).
\]
We call the linearized problem at $\mathbf u^\star$ the system
\begin{equation}
    \label{eq:sistema local linealizado}
    \frac{d \mathbf z}{dt} = \mathbf A \mathbf z \qquad \text{where } \mathbf A = \mathbf J\mathbf f ( \mathbf u^\star ).
\end{equation}
In this type of problems we know that, if $\mathbf f$ is regular enough, and $\mathbf 0$ is a local attractor for \eqref{eq:sistema local linealizado}, then $\mathbf u^\star$ is a local attractor for the nonlinear system \eqref{eq:sistema local no-lin}.
Moreover, the Grobman–Hartman theorem actually tells us that the phase diagrams of \eqref{eq:sistema local no-lin} and \eqref{eq:sistema local linealizado} are locally equivalent in topological terms.

\subsection{Fractional case}
We now study the linear systems
\[
    \dC \mathbf u(t) = \mathbf f(t, \mathbf u(t)), \qquad \dRL \mathbf v(t) = \mathbf g(t, \mathbf v(t)).
\]
Again, existence (and uniqueness) follows from fixed point methods.
We discuss only the linear case, since the linearization results are much more complex (see, for example, \cite{SayevandPichaghchi2015SuccessiveApproximationSurvey}).

\paragraph{Linear equations.}

\begin{figure}[ht!]
    \begin{subfigure}{.5\textwidth}
        \centering
        \includegraphics[]{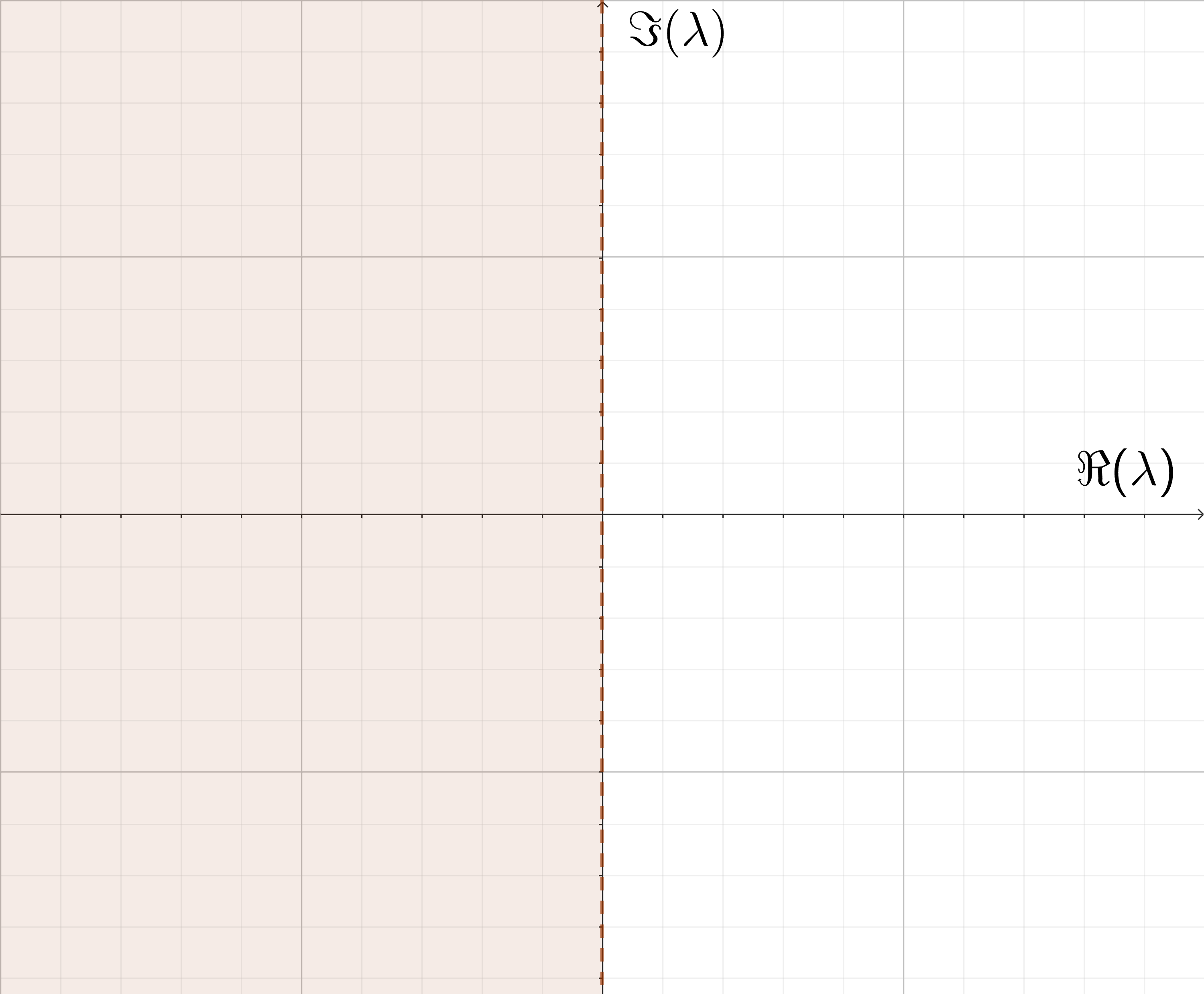}
        \caption{Local case ($\alpha=1$).}
    \end{subfigure}
    \begin{subfigure}{.5\textwidth}
        \centering
        \includegraphics[]{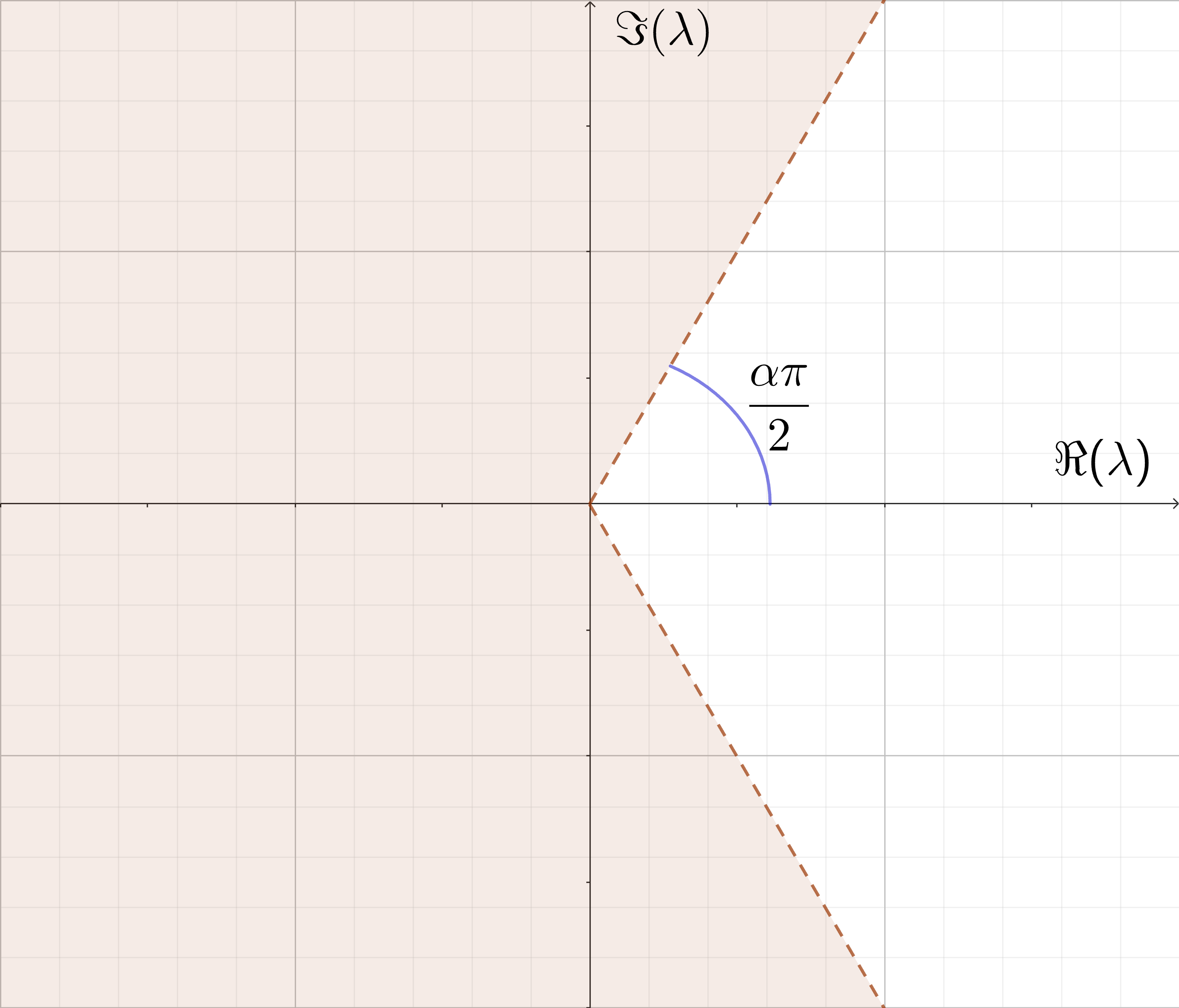}
        \caption{$\dC$ case with $\alpha\in(0,1)$.}
    \end{subfigure}
    \caption{Stability regions (shaded) of autonomous systems.}
    \label{fig:estabilidad}
\end{figure}

Again, take $\mathbf f(t,\mathbf u) = \mathbf g(t, \mathbf u) = \mathbf A \mathbf u$.
As in the local case, if $\mathbf w$ is an eigenvector of $\mathbf A$ with eigenvalue $\lambda$ we arrive at the linear equations $\dC c(t) = \lambda c(t)$ and $\dRL c(t) = \lambda c(t)$ and, therefore, there are particular solutions of the form
\[
    \mathbf u(t) = E_{\alpha} (\lambda t^\alpha) \mathbf w,
    \qquad
    \mathbf v(t) = t^{\alpha - 1} E_{\alpha, \alpha} (\lambda t^\alpha) \mathbf w,
\]
where $E_{\alpha,\beta}$ are the Mittag-Leffler functions defined in \eqref{eq:ML}.
To understand the stability of $\mathbf 0$, we must understand
under what conditions on $\lambda$ we have $t^{\beta - 1} E_{\alpha, \beta} (\lambda t^\alpha) \to 0$ as $t \to \infty$ for $\beta=1$ and $\beta=\alpha$ respectively. The asymptotic behavior of these functions has been studied for years (see \cite{Bateman1985HigherTranscendentalFunctions}).
The necessary and sufficient condition for stability in the $\dC$ case is that for every eigenvalue $\lambda$ of $\mathbf A$ we have $|\arg \lambda| > \frac{\alpha \pi}{2}$ (see \cite{Matignon1996StabilityResultsFractional}).
The necessary and sufficient condition for $\dRL$ is also known and can be found in
\cite{QianLiAgarwalWong2010StabilityAnalysisFractional}. Finally, it is worth noting that fractional problems are “more stable” since even for some eigenvalues with positive real part the problem remains stable. This can be clearly observed in Figure \ref{fig:estabilidad}.

\section{Two applications}
\label{sec:apli}

In this last section we present two situations in which fractional derivatives turn out to be useful tools. To avoid technicalities, the presentation will be at a bird’s-eye view (particularly in the case of the application to finance).

\subsection{The Tautochrone and fractional derivatives}
\label{sec:tautocrona}
N. H. Abel was interested in computing the time that a point mass $m$ subject to the gravitational force would take to slide down a curve in the $x$–$y$ plane described by a certain function $\psi$ via $x = \psi(y)$.
Let us write it parametrically as a curve $\eta \mapsto (\psi(\eta), \eta)$.
Suppose that our point particle of mass $m$ falls under the effect of gravity, of constant value $g$, “stuck” to this curve and starting from rest. To describe its position it suffices to say that at time $t \ge 0$ it is at a height $\phi(t)$. For reasons that will become apparent later, we will use the notation $y = \phi(0)$ and the position vector
\begin{equation*}
    \mathbf z(t) = \Big(\psi(\phi(t)), \phi(t)\Big).
\end{equation*}
Recall the notions of potential energy and kinetic energy for the curve:

\begin{align*}
    E_p(t) & = mg\phi(t),
    \\
    E_c(t) & = \frac 1 2 m \left \|\mathbf v(t)\right \|^2,
\end{align*}
where $\mathbf v(t) = \frac{d \mathbf z}{dt} (t)$ and $\left\|\mathbf v(t)\right \|$ is the speed (whose expression we will compute in detail below).
Moreover, we have the energy conservation formula (neglecting friction)
\begin{align}\label{eq:consener}
    E_p(t)+E_c(t)=E_p(0)+E_c(0).
\end{align}

We are interested in computing the time $T=T(y)$ that the object takes to reach height $0$ (that is, $\phi(T(y))=0$).
We will assume, to simplify the calculations, that the curve is such that $\phi$ is decreasing (i.e. the particle does not go back up).

Since the particle starts from rest, $\phi'(0) = 0$, hence $\mathbf v(0) = \mathbf 0$ and therefore $E_c(0)=0$.
On the other hand, $E_p(0)=mgy$.
It is convenient to write $\mathbf v$ in terms of the arc length defined by the curve described by $\psi$
\[
    s(y) \coloneqq \int_0^y \sqrt{ 1+ \psi'(\eta)^2 }d \eta.
\]
It is immediate to check that
\[
    \| \mathbf v(t) \|^2 = \left\| \Big(\psi'(\phi(t)),1 \Big) \phi'(t) \right\|^2 = \phi'(t)^2 (1 + \psi'(\phi(t))^2) = \phi'(t)^2 s'(\phi(t))^2.
\]
Solving for $\phi'(t)$ in \eqref{eq:consener} we obtain the relation
\[
    -\frac{s'(\phi(t))}{\sqrt{2g(y-\phi(t))}}\phi'(t)=1,
\]
where we have chosen the negative branch of the square root since $s(\phi(t))$ is decreasing.
Integrating both sides between the initial time $t=0$ and the final time $t=T(y)$, and changing variables $\eta=\phi(t)$, we obtain the identity
\begin{equation}
    \label{eq:tautocrona}
    \begin{aligned}
        T(y) & = -\int_0^{T(y)}\frac{s'(\phi(t))}{\sqrt{2g(y-\phi(t))}}\phi'(t)dt    \\
             & =-\int_{\phi(0)}^{\phi(T(y))}\frac{s'(\eta)}{\sqrt{2g(y-\eta)}}d \eta \\
             & =\int_0^{y}\frac{s'(\eta)}{\sqrt{2g(y-\eta)}}d \eta.
    \end{aligned}
\end{equation}
This is known as Abel’s integral formula. In modern notation, we can write this formula in terms of the Caputo fractional derivative $\dC[1/2]$. More precisely, if what is prescribed is a fall–time function from height $y$, $T=T(y)$, one can compute the arc length of the curve (and consequently the curve itself) that yields those times by solving the fractional initial value problem
\begin{align*}
    \dC[1/2]s(y) & = T(y) \sqrt{\frac{2g}{\pi}} , \quad \text{for } y>0, \\
    s(0)         & =0,
\end{align*}
where we have used that $\Gamma(1/2) = \sqrt{\pi}$.
The \textit{Tautochrone problem}
(or isochrone) focused on computing
the curve that would produce a time
of arrival to the ground
independent of the point of the curve at which the object started its descent, that is, $T(y)=k\in \mathbb{R}_+$. It is standard to check, using the results on fractional derivatives acting on monomials (see Section \ref{sec:primera def}), that the solution of
\begin{align*}
    \dC[1/2] [s](y) & = k \sqrt{\frac{2g}{\pi}} , \quad y>0 \\
    s(0)            & =0,
\end{align*}
is given by
\[
    s(y)= \frac{2 k\sqrt{2g} }{\pi}y^{{1}/{2}}.
\]
Since $s'(y)=(1+\psi'(y)^2)^{\frac{1}{2}}$, we can solve for $\psi'$ and use $\psi(0) = 0$ to deduce that
\[
    \psi(y)=\int_0^y\sqrt{\frac{\frac{2g}{\pi^2}k^2}{\eta}-1} \, d \eta.
\]
The curve $(\psi(y),y)$ is precisely a cycloid.
\begin{remark}
    The cycloid also solves other famous problems, for example the brachistochrone problem.
\end{remark}

\begin{remark}
    The usual solution to the tautochrone problem is to set $T(y) = k$ in \eqref{eq:tautocrona}, apply the Laplace transform, and solve.
    This connects with what was explained in Sections
    \ref{sec:EDO lineal no homogenea por Laplace} and \ref{sec:no lineal. nocion de solucion}.
\end{remark}
\subsection{Applications to finance: the Rough Heston model}
\label{sec:finanzas}

The famous Black–Scholes model (1973) attempts to predict the price of financial assets and allows the valuation of call and put options.
These types of models are stochastic, and use the so-called Itô calculus.
An excellent introduction to this theory, to stochastic equations, and to the Black–Scholes model is \cite{Evans2013}.
We now give a very informal presentation of the Black–Scholes model and subsequently of a generalization that uses fractional derivatives.

Let $S_t$ denote the price of a financial asset (for example a stock) at time $t$.
The price $S_0$ is known, but the prices $S_t$ are predicted stochastically.
In more or less general terms, one calls a \emph{stochastic process} a curve $t \in [0,\infty) \mapsto S_t$ where $S_t$ is a random variable.
The deterministic linear growth models described above can be approximated by
\begin{equation*}
    \frac{u(t + \tau) - u(t)}{\tau} \approx \lambda u(t).
\end{equation*}
The idea behind the Black–Scholes model is to try to predict the value at time $t + \tau$ (with $\tau > 0$ small) from $S_t$ by
\[
    S_{t+\tau} - S_t \approx S_t(\tau \mu + \sigma \sqrt{\tau} X_t),
\]
where $\mu$ and $\sigma$ are constants to be calibrated to market data and $X_t$ is a normal distribution $N(0,1)$ such that $X_t$ and $X_s$ are independent if $t \ne s$.
The reason for using $\sqrt{\tau}$ is delicate (it has to do with the central limit theorem), and we recommend \cite{Evans2013} to the interested reader.
We can express the above as a sum
\begin{equation}
    \label{eq:finanzas S discreto}
    S_t \approx S_0 + \mu \tau \sum_{k=0}^{n-1} S_{k \tau} + \sigma \sqrt{\tau} \sum_{k=0}^{n-1} S_{k \tau} X_{k \tau}.
\end{equation}
Define a second stochastic process that will be useful.
The stochastic process $W^{(\tau)}$ described by
\[
    W^{(\tau)}_{t+\tau} = W^{(\tau)}_t + \sqrt{\tau} X_t
\]
is known as a random “walk”.
The limit when $\tau \to 0^+$ of this process is the so-called \textit{Brownian motion}, denoted $W_t$.
The Itô integral allows us to justify “integration with respect to Brownian motion” (see \cite{Evans2013}) as a limit of Riemann sums in \eqref{eq:finanzas S discreto}, which leads us to propose the model
\[
    S_t = S_0 + \mu \int_0^t S_s ds + \sigma \int_0^t S_s d W_s.
\]
For simplicity, it is common to write this in differential notation
\begin{equation}
    \label{eq:Black-Scholes}
    dS_t = \mu S_t dt + \sigma S_t d W_t.
\end{equation}

This model has some limitations. In particular, it was observed that the fact that $\sigma$ (which models the asset’s “volatility”) is constant is a strong limitation for fitting real data.
Many models have been proposed to address this problem. We highlight, for their interest in the fractional field, the stochastic volatility models. The idea is to replace in \eqref{eq:Black-Scholes} the constant $\sigma$ by $\sigma_t = \sqrt{V_t}$, where $V_t$ is a new stochastic process, and write an equation for $V_t$.

A prominent model in this direction was proposed by Heston in 1993.
It consists of taking a second Brownian motion, denoted $B_t$, which can be correlated with $W_t$, and modeling the evolution of $V_t$ with a mean-reverting process with noise, $d V_t = \kappa( \theta - V_t ) dt + \kappa \sqrt{V_t} d B_t$,
or, in integral terms,
$$
    V_t = V_0 + \int_0^t \kappa( \theta - V_s ) ds + \int_0^t \kappa \sqrt{V_t} d B_s.
$$
In the 2010s it was observed that this model also has limitations, and $V_t$ should be “less regular” than what results from this equation
(see \cite{gatheral2014VolatilityRough}).
One of the proposals in this direction
is the \emph{Rough Heston} model, which has had great impact among the financial community in the last 10 years. It consists of replacing the conventional integrals by Riemann–Liouville integrals
$$
    V_t = V_0 + \int_0^t \kappa( \theta - V_s ) K_\alpha (t-s) ds + \int_0^t \kappa \sqrt{V_t}K_\alpha (t-s) d B_s.
$$
One of the arguments that gives validity to this model passes through the modeling of the microstructure of high-frequency trading using Hawkes processes \cite{eleuchMicrostructuralFoundationsLeverage2018a}.
Note that this formula is the stochastic version of the formulation for the ODE associated to $\dC$, that is \eqref{eq:edo Caputo form integral}. We recommend in this direction \cite{eleuch2018QuantitativeFinanceRough} and its references.

For the computation of different quantities associated with this model fractional differential equations arise, which must be solved numerically with some of the methods described above. In particular, the quantitative finance community favors a predictor–corrector method based on Adams methods (see \cite{eleuch2018QuantitativeFinanceRough}).

\section*{Acknowledgments}
The research of Félix del Teso is funded by the Ramón y Cajal contract reference RYC2020-029589-I and the research projects PID2021-127105NB-I00, CNS2024-154515 and CEX2019-000904-S of the AEI, Government of Spain.
The research of David Gómez Castro is funded by the Ramón y Cajal contract reference RYC2022-037317-I and the research project PID2023-151120NA-I0 of the AEI, Government of Spain.


\phantomsection
\end{document}